\documentclass{amsart}
\usepackage{fullpage}
\usepackage{graphicx}
\usepackage{enumerate}
\input pictex.sty

\newtheorem{thm}{Theorem}[section]

\newtheorem{lem}[thm]{Lemma}

\theoremstyle{definition}

\theoremstyle{remark}

\numberwithin{equation}{section}

\newcommand{\To}{\longrightarrow}

\begin{document}

\title{An Algorithm for Map Enumeration}%
\author{Virgil U. Pierce}%
\address{Dept. of Math, The Ohio State University}%
\email{vpierce@math.ohio-state.edu}

\thanks{This work was 
supported in part by NSF grants DMS-0073087,
DMS-0412310, and DMS-0135308.}%

\begin{abstract}
  Bauer and Itzykson showed that associated to each labeled map
  embedded on an oriented Riemann surface there was a group generated
  by a pair of permutations.  From this result an algorithm may be
  constructed for enumerating labeled maps, and this construction is
  easily augmented to bin the numbers by the genus of the surface the
  map is embedded in.   The results agree with the calculations of
  Harer and Zagier of 1-vertex maps; with those of
  Bessis, Itzykson, and Zuber of 4-valent maps; and with
  those of Ercolani, McLaughlin, and Pierce for
  $2\nu$-valent maps.  

  We then modify this algorithm to one which counts unoriented maps or
  Mobius graphs.
  The results in this case agree with the calculation of Goulden and Jackson
   on 1-vertex unoriented maps.
\end{abstract}
\maketitle

\section{Introduction}

The recent work of Ercolani, McLaughlin, and Pierce \cite{EMP} 
 analyzed the fine structure of families of probability
measures on the space of $N \times N$ Hermitian matrices.  
The
principle measures analyzed in \cite{EMP} are of the form
\begin{equation} \label{mu} 
d\mu_t = \frac{1}{Z_N} \exp\left( - N \mbox{Tr}\left[
  V_{t}(M) \right] \right) dM,
\end{equation}
where $dM$ is the product measure
\begin{equation*}
dM = \left( \prod_{1\leq i < j \leq N} d\mbox{Re}(M_{ij}) d\mbox{Im}(M_{ij})
\right) \left( \prod_{1 \leq i \leq N} dM_{ii} \right),
\end{equation*}
the potential is
\begin{equation*}
V_{t}(\lambda) = \frac{1}{2} \lambda^2 + t \lambda^{2\nu},
\end{equation*}
and where $Z_N$ is the normalization of the measure.
The question of interest in that paper was 
the asymptotic structure of $\log(Z_N)$ for
large values of $N$, and its analytic dependence on $t$.

More generally, the partition function of
random matrices, introduced above as the
 normalization factor of the probability measure
(\ref{mu}), is defined as
\begin{equation} \label{part}
Z_N(t_3, t_4, \dots, t_{2\nu}) = \int_{\mathcal{H}_N} \exp\left( - N \mbox{Tr}\left( \frac{1}{2} M^2
+ t_3 M^3 + t_4 M^4 + \dots + t_{2\nu} M^{2\nu} \right) \right) dM,
\end{equation}
where the integral is taken over the
space of $N \times N$ Hermitian matrices.
The partition function has been used extensively as a model for the partition
function of 2D-Quantum Gravity \cite{Witten}, \cite{FGZ}, \cite{BDiFG}.  

Derivatives of $Z_N(t_3, t_4, \dots, t_{2\nu})$ evaluated at $t_j= 0$
capture the Gaussian moments $\int_{\mathcal{H}_N} \mbox{Tr}(M^j)^k
d\mu_0(M)$.  These moments are evaluated with the Wick lemma
in terms
of pair correlators of the matrix entries.  In 
analogy to Feynman diagrams, there is a correspondence between labeled
maps and the terms in the Wick expansion.
This fact led to the conjecture \cite{BIZ} that the logarithm of the
partition function (\ref{part}) has an asymptotic expansion of the
form 
\begin{equation*}
\log\left[ Z_N(t_3, t_4, \dots, t_{2\nu})/Z_N(0) \right] = N^2 e_0(t_3, t_4,
\dots, t_{2\nu}) + e_1(t_3, t_4, \dots, t_{2\nu}) + \dots
\end{equation*}
where the coefficients $e_g(t_3, t_4, \dots, t_{2\nu})$ are
generating functions for counting labeled, connected, oriented maps
embedded on a Riemann surface of genus $g$.

A map $D$ on a compact, oriented connected surface $X$ is a pair $D=(K(D),
[\iota])$ where 
\begin{enumerate}
\item $K(D)$ is a connected 1-complex;
\item $[\iota]$ is an isotopical class of inclusions $\iota : K(D) \To X$;
\item the complement of $K(D)$ in $X$ is a disjoint union of open cells
  (faces);
\item the complement of $K_0(D)$ (vertices) in $K(D)$ is a disjoint union of
  open segments (edges).
\end{enumerate}
The number of edges connected to a vertex $z \in K_0(D)$ is called the degree
of $z$.  We may represent maps as fat or ribbon graphs embedded on an oriented
connected surface $X$.  See also the theory of dessins d'enfants \cite{Sc94}.

Ercolani and McLaughlin \cite{emcl_2003} rigorously showed that $\log(Z_N)$
possesses an asymptotic expansion of the form 
\begin{equation*}
\frac{1}{N^2} \log\left( \frac{Z_N(t_3, t_4, \dots, t_{2\nu})}{Z_N(0)} \right) 
= e_0(t_3, t_4, \dots, t_{2\nu}) + \frac{1}{N^2} e_1(t_3, t_4, \dots,
t_{2\nu}) + \dots
\end{equation*}
in a non-trivial $(t_3, t_4, \dots, t_{2\nu})$ domain.  The functions 
$e_g(t_3, t_4, \dots, t_{2\nu})$ are analytic in a
neighborhood of $0$ and are counting functions for  labeled maps.

The function 
\begin{equation*}
e_g(t_3, t_4, \dots, t_{2\nu}) = \sum_{j_3, j_4, \dots, j_{2\nu} = 0}^\infty
\kappa_g(j_3, j_4, \dots, j_{2\nu}) t_3^{j_3} t_4^{j_4} \dots t_{2\nu}^{j_{2\nu}},
\end{equation*}
where $\kappa_g(j_3, j_4, \dots, j_{2\nu})$ is the number of
labeled, genus $g$, maps with $K_0(D)$ containing:
$j_3$ vertices of degree $3$, 
$j_4$ vertices of degree $4$, $\dots$, $j_{2\nu}$ vertices of degree
$2\nu$.

Harer and Zagier \cite{Harer_Zagier} solved this enumeration problem
in the special case of 1-vertex maps (or monopoles).  
Let $G_n(N)$ be the generating function for oriented maps with a
single vertex, that is 
\begin{equation*}
G_n(N) = N^{1+n} \sum_{ g \geq 0} a_{n, g} N^{-2g},
\end{equation*}
where $a_{n, g}$ is the number of oriented maps with 1 vertex, and $n$
edges, which are embedded into a genus $g$ surface.
This generating function is a finite series in $N$ and is given by the
matrix integral 
\begin{equation} \label{G_n(N)_int}
G_n(N) = A_N \int_{\mathcal{H}_N} \mbox{Tr}(M^{2n}) \exp\left[
  -\frac{1}{2} \mbox{Tr}(M^2) \right] dM,
\end{equation}
where the integral is taken over the space of $N\times N$ Hermitian
matrices, and  
\begin{equation*}
A_N  = \left[ \int_{\mathcal{H}_N} \exp\left[ -\frac{1}{2}
  \mbox{Tr}(M^2)\right] dM \right]^{-1}.
\end{equation*}

\begin{thm}[Harer-Zagier \cite{Harer_Zagier}]
Harer and Zagier explicitly evaluated (\ref{G_n(N)_int}), to find that
\begin{equation*}
G_n(N) = \frac{(2n)!}{2^n n! } \sum_{k=0}^n 2^k \binom{n}{k} \binom{N}{k+1}.
\end{equation*}
\end{thm}

Bauer and Itzykson \cite{BI} show that the data of a labeled 
map is equivalent to a pair of permutations satisfying some compatibility
conditions. 
This result gives an efficient method to
determine the genus of a map given by its pair of permutations.  
What emerges is an algorithm for
counting the number of genus $g$ maps with a particular structure
on $K_0(D)$.
The idea is to  count the number of
pairs of permutations satisfying the characterization which produce a genus $g$
map.   Zvonkin \cite{Zvonkin2} expanded on this idea.

The particular partition function $Z_N(0, 0, \dots, t_{2\nu})$ has been the 
subject of an extensive calculation \cite{EMP}.  The conclusion of this
calculation is a prescription for calculating $e_g(0, 0, \dots, t_{2\nu})$.
Explicit formulas for these functions are found for low values of $g$ ($g=0,
1, 2, 3$) in terms of an auxiliary function.  In the case of $g=0$ this result
produces a closed form expression for the Taylor coefficients $\kappa^{(2
\nu)}_0(0, 0, \dots, j_{2\nu})$.  In the other cases worked out explicitly this
result is used to find the first few Taylor coefficients of each of the
functions.  

A purely combinatoric method for producing an expression for
$e_0(t_3, t_4, \dots, t_{2\nu})$ in terms of generating functions of
fundamental objects is given in \cite{BDiFG}.

The algorithm detailed in this paper  has been 
developed as an independent method
for checking the 
detailed work on the fine structure of $e_g(0, 0, \dots, t_{2\nu})$ carried
out in \cite{EMP} and of deriving the constants of integration for the
formulas derived in that work, and in the author's doctoral thesis.   
As of the writing of this paper the calculation has been done for the Taylor
coefficients satisfying $\nu j \leq 20$ and we have found exact agreement for
all of these terms.

We will conclude this paper with a modification of our method for counting the
number of unoriented maps or Mobius graphs.  These are maps which are
embedded into an unoriented surface.  To be more precise the objects
we will count are Mobius graphs, ribbon graphs where the edges are
allowed to twist.
This is the setting of interest if the
partition function (\ref{part}) is replaced with
\begin{equation}\label{ZNGOE}
Z^{(1)}_N(t_3, t_4, \dots, t_{2\nu}) = \int_{\mathcal{S}_N} \exp\left[
  -N\mbox{Tr}\left( \frac{1}{4} M^2 + t_3 M^3 + t_4 M^4 + \dots +
  t_{2\nu} M^{2\nu} \right)\right] dM ,
\end{equation}
where the integral is  over $N \times N$ symmetric
matrices.  This function is called the partition function of the
Gaussian orthogonal ensemble (GOE).

In this case the hypothesis is that 
\begin{equation}\label{GOE_expansion}
\frac{1}{N^2} 
\log\left( \frac{Z^{(1)}_N(t_3, t_4, \dots, t_{2\nu})}{Z_N(0)} \right) 
= E_0(t_3, t_4, \dots, t_{2\nu}) + \frac{1}{N} E_1(t_3, t_4, \dots, t_{2\nu}) 
+ \frac{1}{N^2} E_2(t_1, t_2, \dots, t_{2\nu} + \dots
\end{equation}
where (at least formally)
\begin{equation*}
E_\chi(t_3, t_4, \dots, t_{2\nu}) = \sum_{j_3, j_4, \dots, j_{2\nu} =
  0}^\infty \kappa_{\chi}(j_3, j_4, \dots, j_{2\nu}) t_3^{j_3}
  t_4^{j_4} \dots t_{2\nu}^{j_{2\nu}}, 
\end{equation*}
where  $\kappa_\chi(j_3, j_4, \dots, j_{2\nu})$ is the number of
labeled, Euler characteristic $\chi$, unoriented maps with $K_0(D)$
containing: $j_3$ vertices of degree 3, $j_4$ vertices of degree 4,
$\dots$, $j_{2\nu}$ vertices of degree $2\nu$.  
The justification of this enumeration problem is given in
\cite{Goulden_Jackson} and \cite{Mulase}.  

The existence of the expansion (\ref{GOE_expansion}) 
is not rigorously shown to date.  
In
particular a proof is needed of the analyticity of $E_\chi$.  We will
give a heuristic argument for the counting property of $E_\chi$.  

Goulden and Jackson \cite{Goulden_Jackson} generalized the calculation
of Harer and Zagier \cite{Harer_Zagier} for 1-vertex maps (monopoles).  
Let $F_n(N)$ be the generating function for unoriented maps with a
single vertex, that is 
\begin{equation*}
F_n(N) = N^{n-1} \sum_{\chi \leq 2} f_{n, \chi} N^{\chi} ,
\end{equation*}
where $f_{n, \chi}$ is the number of unoriented maps with 1 vertex,
and $n$ edges, which are embedded into a genus $g$ surface.
This generating function is a finite series in $N$ and is given by the
matrix integral
\begin{equation} \label{F_n(N)_int}
F_n(N) = B_N \int_{\mathcal{S}_N} \mbox{Tr}(M^{2n}) \exp\left[
  -\frac{1}{4} \mbox{Tr}(M^2) \right] dM, 
\end{equation}
where the integral is taken over the space of $N\times N$ symmetric
matrices, and 
\begin{equation*}
B_N = \left[ \int_{\mathcal{S}_N} \exp\left[ -\frac{1}{4} \mbox{Tr}(M^2)
  \right] dM \right]^{-1}.
\end{equation*}
  
\begin{thm}[Goulden-Jackson \cite{Goulden_Jackson}] \label{GJ_thm}
Goulden and Jackson explicitly evaluated (\ref{F_n(N)_int}), to find that 
\begin{equation*}
F_n(N) = n! \sum_{k=0}^n 2^{2n-k} \sum_{r=0}^n \binom{n-1/2}{n-r}
\binom{k+r-1}{k} \binom{ (N-1)/2}{r} + \frac{(2n)!}{2^n n!}
\sum_{k=0}^n 2^k \binom{n}{k} \binom{N-1}{k+1}.
\end{equation*}
\end{thm}

Our algorithm agrees with the numbers given by $F_n(N)$ for all cases checked.
What we do not have is a method for generalizing Theorem \ref{GJ_thm}
to larger families of unoriented maps.

In section \ref{wick_sec} we outline the heuristic argument which
motivates using the terms of the asymptotic expansion of the partition
function (\ref{part}) as the generating functions of the number of
maps of genus $g$.  In section \ref{group_sec} we illustrate the
result of \cite{BI} which gives a one-to-one correspondence between
maps and pairs of 
permutations, of a restricted type.   So the 
group structure provides an algorithm for counting the number of maps
of a specified type.  Section \ref{results_sec} is a collection of
tables summarizing some of the calculations we have done with this
algorithm.  It is provided as a means of checking with other known
results such as the 1-face maps \cite{Harer_Zagier}, \cite{Zvonkin}
and those in \cite{BIZ} and \cite{EMP}.  In addition the tables
contain all the numbers which were needed to evaluate the constants of
integration found in \cite{EMP}.
In section \ref{sym_sec} we redo the heuristic argument of section
\ref{wick_sec} in the setting of an integral over the symmetric
matrices rather than the Hermitian matrices.  We show that the
equivalent counting problem we are concerned with is that of counting
maps embedded into an unoriented surface.  In section
\ref{unoriented} we generalize the results of section \ref{group_sec} and
associate to each unoriented map a triple of permutations. 
This gives a one-to-one correspondence between unoriented maps and
triples of permutations of the given type.

\section{\label{wick_sec} The Wick lemma and Hermitian matrix integrals}

We will first outline the heuristic argument that the $e_g(t_3,
t_4, \dots, t_{2\nu})$ are counting functions of labeled genus $g$
maps (see \cite{BIZ} and \cite{Zvonkin}).  
This calculation centers around the Wick lemma for Gaussian
expectation's:
The Wick lemma states 
\begin{lem} \label{WickLemma}
If $\langle \cdot \rangle$ is a Gaussian expectation, and $l_1, l_2,
\dots, l_{2n}$ are homogeneous linear functions then
$$ \langle l_1 l_2 \dots l_{2n} \rangle 
= \sum \langle l_{i_1} l_{j_1} \rangle \langle l_{i_2} l_{j_2} \rangle
\dots \langle l_{i_n} l_{j_n} \rangle,
$$ 
where the sum is taken over the set $1 = i_1< i_2 < \dots < i_n < 2n$ and 
$i_k < j_k$.  
\end{lem}

We are working with the Gaussian Expectation
\begin{equation*}
\langle f(M) \rangle = C_N \int_{\mathcal{H}_N} f(M) \exp\left[ - N
  \frac{1}{2} \mbox{Tr}\left(M^2\right) \right] dM ,
\end{equation*}
where the integral is taken over the space of $N \times N$ Hermitian
matrices, and $C_N$ is a normalizing constant.  The quadratic
expectations of this probability are 
\begin{equation} \label{GUE_exp}
\langle M_{i j} M_{k l} \rangle = \frac{1}{N} \delta_{i l} \delta_{j k} .
\end{equation}

Consider the partition function 
\begin{equation} \label{Hermitian_Matrix}
Z_N(t_4) = C_N \int_{\mathcal{H}_N} \exp\left[ -N \mbox{Tr}\left(
  \frac{1}{2} M^2 + t_4 M^4 \right) \right] dM .
\end{equation}
We may expand out the $t_4$ expression as a Taylor Series and
(formally at least) commute this sum with the integration to write 
\begin{equation}\label{eq1}
Z_N(t_4) = \sum_{n=0}^\infty \frac{ (-t_4)^n}{n!} N^n \langle
\left[ \mbox{Tr}\left(M^4\right)\right]^n \rangle . 
\end{equation}
Write (\ref{eq1}) as a function of the entries of $M$:
\begin{align} \label{eq2}
Z_N(t_4) &= \sum_{n=0}^\infty \frac{(-t_4)^n}{n!} \sum_{i_m, j_m, k_m,
    l_m} \langle
    \prod_{m=1}^n M_{i_m
    j_m} M_{j_m k_m} M_{k_m l_m} M_{l_m i_m} \rangle \\
&= \nonumber
1 + (-t_4) \sum_{i, j, k, l} 
\langle M_{i j} M_{j k} M_{k l} M_{l i} \rangle  + \dots
\end{align}

If we use Wick's lemma on the $n=1$ term of (\ref{eq2}) 
we find that the contribution is given by 
\begin{equation}
 \langle M_{i j} M_{j k} \rangle \langle M_{kl} M_{li} \rangle +
 \langle M_{i j} M_{k l} \rangle \langle M_{j k} M_{l i} \rangle + 
\langle M_{i j} M_{l i} \rangle \langle M_{j k} M_{k l} \rangle
\end{equation}

\begin{center}
\begin{figure}[h]
\scalebox{0.9}{\includegraphics{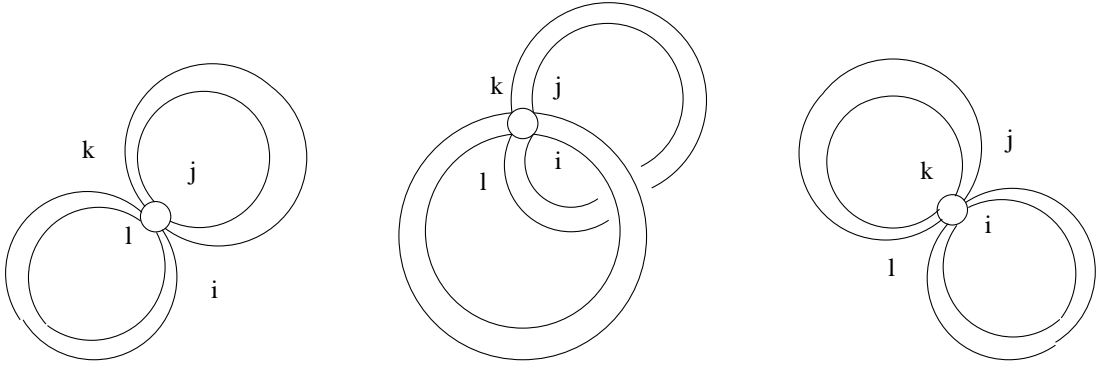}}
\caption{\label{fig1} The three ribbon graphs representing the maps
  associated to the $n=1$ term of (\ref{eq2}).  }  
\end{figure}\end{center}

Each of these terms is non-zero if $i, j, k,$ and $l$ satisfy the
conditions for the corresponding quadratic expectations to be
non-zero.  In this example we find the conditions are respectively 
$ i = k$, $ i = l = k = j$, and $j = l$.  In an analogy to Feynman
diagrams these reductions are
cataloged by families of maps (see for example Figure \ref{fig1}). 
Moreover the number of relations between the indices is the number of
faces of the corresponding map, while the contribution to the
expectation $\langle \exp\left[-N\mbox{Tr}\left(M^4\right)\right]
\rangle$ is given by the Euler Characteristic of the map.  

 The fact that the maps organizing the contributions to the Hermitian
 matrix integral (\ref{Hermitian_Matrix}) are oriented
is because of (\ref{GUE_exp}).  

The final step is to take a logarithm of $Z_N(t_4)$ to restrict our
consideration to connected maps.  This argument easily generalizes to
include other times.  

\section{\label{group_sec} Group data associated to each map $[D]$}

We will now outline the group data associated to a map, the data comes
in the form of a pair of permutations.  These two permutations
provide a mechanism for counting the number of maps by computing the
number of permutations satisfying some conditions.  Conveniently this
permutation structure also provides an efficient method of computing
the structure of the faces of the map and therefore of calculating the
Euler Characteristic of the map.

  Let $E$ be the number of edges of $K(D)$.  
We will view the permutations of $\mathcal{S}_{2E}$ as acting on the set of
darts of $K(D)$, which is represented by 
\begin{equation*}
\Omega = \Omega_1 \cup \Omega_2, 
\end{equation*}
where 
\begin{equation*}
\Omega_1 = \left\{ (v, e): 
\mbox{ $v$ is a vertex and $e$ is an edge with two distinct vertices
  one of them being $v$ } 
\right\}, 
\end{equation*}
and 
\begin{equation*}
\Omega_2 = \left\{ (v, e, \pm): 
\mbox{ $v$ is a vertex and $e$ is an edge with a single vertex $v$} \right\}.
\end{equation*}
The element $(v, e) \in \Omega_1$ represents the dart based at $v$ and
going along $e$.  The element $(v, e, \pm) \in \Omega_2$ represents
the dart based at $v$ going along $e$ in the counterclockwise
(resp. clockwise) orientation.  
For each edge there are two darts therefore $\big| \Omega \big| =  2
E$ and 
we can think of $\mathcal{S}_{2E}$ as acting by permutations on the
set $\Omega$.

Given a map, $(K(D),[\iota])$, we
define a subgroup of $\mathcal{S}_{2E}$ generated by two permutations
$\langle \sigma, \tau \rangle$. 
The orientation on $X$ induces (via $[\iota]$) a cyclic ordering on
the darts attached to each vertex; the first permutation $\sigma$ is
given by this action.  Explicitly, $\sigma$ maps the element $(v, e)
\in \Omega_1$ to the element $(v, \tilde{e}) \in \Omega_1$ or $(v,
\tilde{e}, \pm) \in \Omega_2$ where $\tilde{e}$ is the
edge counter clockwise in the orientation at $v$ from $e$. 
Likewise, $\sigma$ maps the element $(v, e, \pm) \in \Omega_2$ to the
element $(v, \tilde{e}) \in \Omega_1$ or $(v, \tilde{e}, \pm) \in
\Omega_2$ where $\tilde{e}$ is the edge counter clockwise in the
orientation at $v$ from $e$.  
  The second permutation $\tau$, is given explicitly as 
the permutation which acts on $\Omega_1$ by sending $(v, e)$ to
$(\tilde{v}, e)$ where  $\tilde{v}$ is the other endpoint of $e$; and
$\tau$ acts on $\Omega_2$ by sending $(v, e, \pm)$ to $(v, e, \mp)$.

One sees that $\sigma$ is a product of independent cycles;
the length of the cycles in $\sigma$ correspond to the degrees of the
vertices in $K_0(D)$.  
 The
permutation $\tau$ is a product of $E$ disjoint 2-cycles.  

The orbit of $< \sigma, \tau >$ (the group generated by
$\sigma$ and $\tau$) of any element in
$\Omega$, is all of $\Omega$ because $K(D)$ is connected.  

To compute the genus of a map from $(\sigma, \tau)$ we compute the
genus of the associated graph from the Euler Characteristic
\begin{equation}\label{euler}
\chi = V - E + F;
\end{equation}
where $V$, the number of vertices, is the number of cycles in $\sigma$;
where $E$, the number of edges, is $1/2$ the number of elements acted upon by
$\sigma$ and $\tau$, or is the number of 2-cycles in $\tau$; and where
$F$, the number of faces, is the number of cycles of the permutation
 $\sigma \circ \tau$ (this is merely computing each face by finding all the
edges which border it).  

This leads to the theorem:
\begin{thm}[Bauer and Itzykson \cite{BI}] 
There is a one-to-one correspondence between connected maps $D = (K(D),
[\iota])$ and pairs of permutations ($\sigma$, $\tau$), where $\sigma$
is a product of disjoint cycles, and $\tau$ a product of $E$ disjoint
2-cycles, which satisfy connectedness.  
\end{thm}

To complete the proof we show that given ($\sigma$, $\tau$) as above
one may determine $D$.
The condition of connectedness is that the orbit of the group $\langle \sigma,
\tau \rangle \cdot 1$ is all $2E$ letters.
The
following steps produce a map:
\begin{enumerate}
\item Let $K_0(D)$ be the set of cycles of $\sigma$;
\item Let $K_1(D)$ be the set of 2-cycles of $\tau$;
\item Each 2-cycle in $\tau$ ties together two points of $K_0(D)$
  identifying the two endpoints of that element of $K_1(D)$;
\item The Euler Characteristic (\ref{euler}) gives the genus of $X$;
\item Each cycle in $\sigma$ induces the orientation on $X$
 in a neighborhood of $\iota(K_0(D))$, this is sufficient to determine both
 $X$ and $\iota$ up to isotopy (see section \ref{unoriented}).  
\end{enumerate}

The punchline is that what we now have is an algorithm for computing
$\kappa_g(j_3, j_4, \dots, j_{2\nu})$ for finite $g$, $j_3$, $j_4$, $\dots$, $j_{2\nu}$:
Let 
\begin{equation*}
2 E = 3 j_3 + 4 j_4 + \dots + 2\nu j_{2\nu}.
\end{equation*}
Fix $\sigma$ to be a permutation formed by $j_3$ cycles of length 3, $j_4$
cycles of length 4, $\dots$, $j_{2\nu}$ cycles of length $2\nu$ formed from a
permutation of $n$ letters.  Then we
choose each product of $E$ disjoint
2-cycles in $\mathcal{S}_{2E}$.    
Check if $(\sigma, \tau)$ is connected (by verifying that the orbit of
$<\sigma, \tau> \cdot 1$ is all $2E$ letters).  If $(\sigma, \tau)$ is
connected compute the genus as above (\ref{euler}) 
with $\chi = 2 - 2g$ and bin the result.

\subsection{Example \label{example_section}}

As an example we will show how this algorithm works by computing 
$\kappa_g( j_3 = 2, j_4 = 1)$.  
Fix $\sigma = (1 \; 2 \; 3)(4 \;5\; 6)(7 \;8 \;9 \;10)$.  Then choose each
$\tau$ which is a 
product of 2-cycles of $\{ 1, 2, \dots, 10 \}$.  

For example, we might choose $\tau = (1\; 2 )(3 \;4 ) (5 \;6)( 7 \;8 ) (9 \;10)$; 
for this choice of $\tau$ the orbit of $<\sigma, \tau>$ forms two disjoint
sets $\{ 1, 2, \dots, 6 \}$ and $\{ 7, 8, 9, 10\}$.  

Choose $\tau = (1 \;2)(3\; 4)(5\; 8) (6 \;7) (9 \;10)$; for this choice of
$\tau$ the 
orbit of $<\sigma, \tau>$ is the set $\{ 1, 2, \dots, 10\}$, therefore
$(\sigma, \tau)$ represent a map.  We compute the genus of this map by
computing the number of distinct orbits of $\sigma \cdot \tau$.  We find that 
$\sigma \cdot \tau = (1 \;3 \;5\; 9\; 7\; 4) (2) (6 \;8) (10)$, therefore the
Euler 
characteristic is $3 - 5 + 4 = 2$ and the genus is $0$.  

Choose $\tau = (1 \;2) (3 \;4) (5 \;7) (6 \;9) (8 \;10)$; for this choice of
$\tau$ the 
orbit of $<\sigma, \tau>$ is the set $\{ 1, 2, \dots, 10\}$, therefore
$(\sigma, \tau)$ represent a map.  We find that $\sigma \cdot \tau = (1 \;3 \;5 \;8
\;7\;6 \;10\; 9\; 4) (2)$, therefore the Euler characteristic is $3 - 5 + 2 = 0$ and the
genus is $1$.

\section{Results \label{results_sec}  }

In tables \ref{table1}-\ref{table3}
 we present a selection of the results of this calculation, included are
the numbers which were needed for the calculations in \cite{EMP}.  

The algorithm can also be used to compute the number of maps of mixed type.
For the example from section \ref{example_section} of 2 vertices of degree 3
and one of degree 4 we find 432 genus 0 maps and 468 genus 1 maps.    
For the number of maps with one vertex of degree 3, one of degree 4,
and one of degree 5: 2160 (genus 0), 6480 (genus 1), and 1440 (genus 2).

\begin{table}[h]
\caption{ \label{table1}
First the number of one vertex maps found by \cite{Harer_Zagier}
 (whose duals are the well studied
1-face maps \cite{Zvonkin}) }
\begin{tabular}[b]{|r||r|l|l|l|l|l|l|}
\hline
Degree & Genus         & 0     & 1       & 2        & 3         & 4         & 5 \\ \hline\hline
4      &               & 2     & 1       &          &           &           &   \\
6      &               & 5     & 10      &          &           &           &   \\
8      &               & 14    & 70      & 21       &           &           &   \\
10     &               & 42    & 420     & 483      &           &           &   \\
12     &               & 132   & 2310    & 6468     & 1485      &           &   \\
14     &               & 429   & 12012   & 66066    & 56628     &           &   \\
16     &               & 1430  & 60060   & 570570   & 1169740   & 225225    &   \\ 
18     &               & 4862  & 291720  & 4390386  & 17454580  & 12317877  & \\
20     &               & 16796 & 1385670 & 31039008 & 211083730 & 351683046 &59520825
\\\hline
\end{tabular}
\end{table}

\begin{table}[h]
\caption{\label{table2} The number of two vertex maps}
\begin{tabular}[b]{|r||r|l|l|l|l|l|}
\hline
Degree & Genus & 0      & 1       & 2        & 3         & 4 \\ \hline\hline
3      &       & 12     & 3       &          &           &   \\
4      &       & 36     & 60      &          &           &   \\
5      &       & 180    & 600     & 165      &           &   \\
6      &       & 600    & 4800    & 4770     &           &   \\
7      &       & 2800   & 34300   & 81340    & 16695     &   \\
8      &       & 9800   & 215600  & 1009400  & 781200    &   \\
9      &       & 44100  & 1323000 & 10478160 & 19158300  & 3455865 \\
10     &       & 158760 & 7408800 & 94091760 & 333774000 & 218402730 \\
\hline
\end{tabular}
\end{table}

\begin{table}[h]
\caption{\label{table3}
For the case of degree 4 maps we sum up our calculations in the following
table (these numbers agree with those found in \cite{BIZ} and
\cite{EMP})
}
\begin{tabular}[b]{|r||r|l|l|l|l|}
\hline
Vertices & Genus & 0        & 1         & 2         & 3 \\ \hline\hline
1        &       & 2        & 1         &           &   \\
2        &       & 36       & 60        &           &   \\
3        &       & 1728     & 6336      & 1440      &   \\ 
4        &       & 145152   & 964224    & 770688    &   \\
5        &       & 17915904 & 192098304 & 348033024 & 58060800 \\
\hline
\end{tabular}
\end{table}

\section{\label{sym_sec} Wick lemma and symmetric matrix integrals}

This section is motivated by 
the partition function $Z_N^{(1)}$ given by (\ref{ZNGOE}).  

We will now show, using the Wick lemma, that the asymptotic expansion
of $\log\left( Z^{(1)}_N \right)$ enumerates unoriented maps.  This will be a
formal calculation only, the necessary analyticity of the terms of this 
expansion has not been rigorously shown.  We present the calculation
 here as
a motivation for the counting problem we are computing.

We are working with the Gaussian Expectation 
\begin{equation*}
\langle f(M) \rangle = C_N \int_{\mathcal{S}_N} f(M) \exp\left[ -N
  \frac{1}{4} \mbox{Tr}(M^2)\right] dM,
\end{equation*}
where the integral is taken over the space of $N \times N$ symmetric
matrices, and $C_N$ is a normalizing constant.  
The quadratic expectations of this probability are 
\begin{equation} \label{quad}
\langle M_{ij} M_{kl} \rangle = \frac{1}{N} \left[ \delta_{il}
  \delta_{jk} + \delta_{ik} 
\delta_{jl} \right].
\end{equation}

Consider the partition function 
\begin{equation*}
Z^{(1)}_N(t_4) = C_N \int_{\mathcal{S}_N} \exp\left[ - N \mbox{Tr}\left(
  \frac{1}{4} M^2 + t_4 M^4 \right) \right] dM.
\end{equation*}
We may expand out the $t_4$ expression as a Taylor Series and
(formally at least) commute this sum with the integration to write
\begin{equation*}
Z^{(1)}_N(t_4) = 
\sum_{n=0}^\infty \frac{(-t_4)^n}{ n!} N^n \langle \mbox{Tr}(M^4)^n
\rangle .
\end{equation*}
Write this expression as a function of the entries of $M$: 
\begin{align*}
Z^{(1)}_N(t_4) &= \sum_{n=0}^\infty \frac{(-t_4)^n}{n!} N^n \sum_{i_m, j_m,
  k_m, l_m} \langle 
\prod_{m=1}^n 
M_{i_m j_m}
M_{j_m k_m} M_{k_m l_m} M_{l_m i_m} \rangle \\ 
&= 1 + (-t_4) \sum_{i, j, k, l} \langle M_{i, j} M_{j, k} M_{k, l}
  M_{l, i} \rangle + \dots
\end{align*}

If we use Wick's Lemma on the $n=1$ term we find that the contribution
is given by 
\begin{equation}  \label{Wick}
\langle M_{i j} M_{j k} \rangle \langle M_{k l} M_{l i} \rangle 
+ \langle M_{i j} M_{k l} \rangle \langle M_{j k} M_{l i} \rangle
+ \langle M_{i j} M_{k i} \rangle \langle M_{j k} M_{k l} \rangle
\end{equation}
Each of these terms is non-zero if $i, j, k$ and $l$ satisfy the
conditions for the corresponding quadratic expectations to be
non-zero.  As in the case of oriented maps these reductions may be
cataloged by families of maps.  In this case the maps are allowed to
have edges that reverse orientation because of formula (\ref{quad}).

\begin{center}
\begin{figure}[h] 
\scalebox{0.8}{\includegraphics{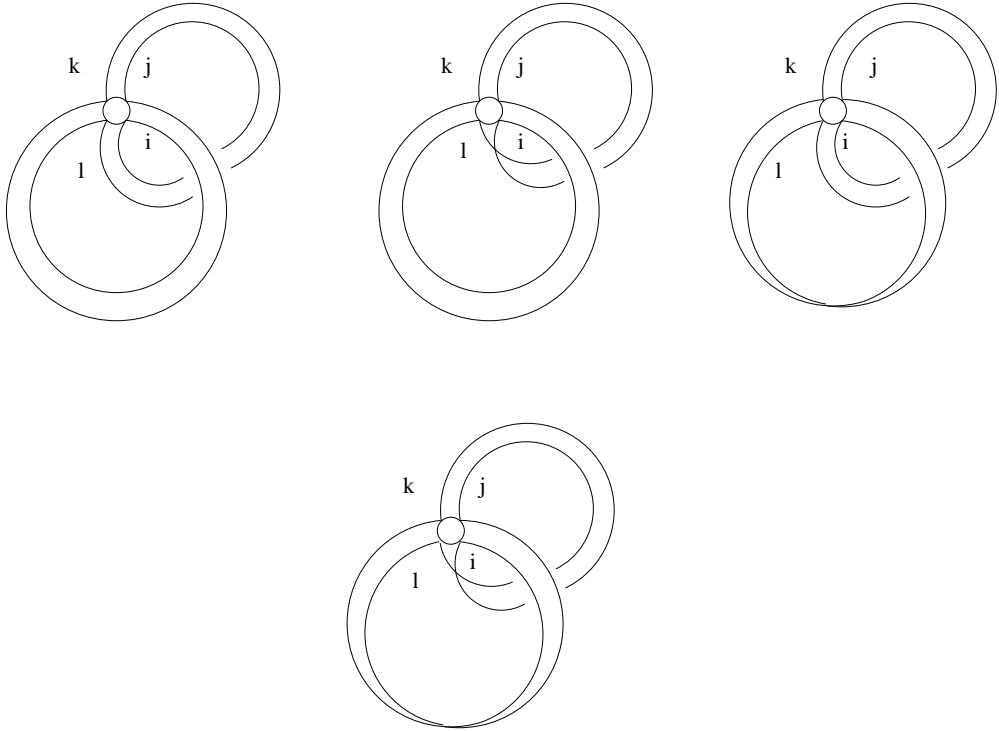}}
\caption{\label{fig2} Four Mobius graphs corresponding to the unoriented maps whose
  underlying oriented map is the single, genus one, oriented map with one
  vertex of degree four.}
\end{figure}\end{center}

Focus on the middle term of (\ref{Wick})
\begin{equation}
\langle M_{i j} M_{k l} \rangle \langle M_{j k} M_{l i} \rangle.
\end{equation}
There are four cases which give that this contribution is non-zero:
$ i = l = k =j$, $i = k = l = j$, $i = l = k = j$, or $(i = k, j=l)$,
these correspond to the four unoriented maps represented by the Mobius
graphs in Figure \ref{fig2}.  

The final step is to take a logarithm of $Z^{(1)}_N(t_4)$ to restrict to
connected unoriented maps.

\section{Unoriented maps \label{unoriented} }

We will lift the algorithm for counting labeled connected oriented
maps to count labeled connected unoriented maps.  We will first go
over the modifications which must be made to the group data to account
for the independence that the local orientations of each vertex have.

An unoriented map on a compact unoriented connected surface $X$ 
is a pair $D = ( K(D), [\iota]) $ where 
\begin{enumerate}
\item $K(D)$ is a connected 1-complex;

\item $[\iota]$ is an isotopical class of inclusions $\iota: K(D) \to X$;

\item the complement of $K(D)$ in $X$ is a disjoint union of open
  cells (faces);

\item the complement of $K_0(D)$ (vertices) in $K(D)$ is a disjoint
  union of open segments (edges).

\end{enumerate}
In the same way that ribbon graphs are associated to an oriented map
we associate a Mobius graph (ribbon graph with orientation reversing
and orientation preserving edges) to an unoriented map.

Given an unoriented map $D$, there are three permutations associated
to the data contained in $D$. 
To begin, for each vertex, $v \in K_0(D)$  choose a local
orientation, $\vec{n}_v$ induced by the embedding $\iota$.
 We will view the four permutations as
acting on the set of darts $\times \mathbb{Z}_2$, given explicitly as: 
$
\Omega = \Omega_1 \cup \Omega_2,
$
where 
\begin{equation*}
\Omega_1 = \left\{ (v, e, \pm \vec{n}_v):  \mbox{$v$ is a vertex and
  $e$ is an edge with two distinct vertices one of them being $v$}
\right\},
\end{equation*}
and
\begin{equation*}
\Omega_2 = \left\{ (v, e, \pm \vec{n}_v, \pm): \mbox{$v$ is a vertex
  and $e$ is an edge with a single vertex $v$} \right\}.
\end{equation*}
The element $(v, e, \pm \vec{n}_v) \in \Omega_1$ 
represents the dart based
at $v$ pointing along $e$ together with the local orientation $\vec{n}_v$
(resp. $-\vec{n}_v$).  The element $(v, e, \vec{n}_v; \pm) \in \Omega_2$
represents the dart based at $v$ which traverses $e$ in the counter
clockwise direction with respect to the local orientation 
(resp. clockwise direction with respect to the local orientation)
together with the local orientation $\vec{n}_v$.   Likewise the element
$(v, e, -\vec{n}_v; \pm) \in \Omega_2$ represents the dart based at $v$ which
traverses $e$ in the counter clockwise direction with respect to the
local orientation (resp. clockwise direction with respect to the local
orientation) together with the normal vector $-\vec{n}_v$.  

  An edge is said to be
orientation reversing if $\vec{n}_v$ at $v$ on
the edge is $-\vec{n}_{\tilde{v}}$ 
after translating to $\tilde{v}$ along the edge; 
otherwise we say that the edge is orientation
preserving.  This is a property of both the map and the choice of
local orientations.

We find three permutations $(\phi, \sigma, \tau)$ :

On $\Omega_1$ we define the action
of $(\phi,
\sigma, \tau)$ by:
\begin{itemize}
\item
$\phi$ acts by sending $(v, e, \pm \vec{n}_v)$ to
$(v, e, \mp \vec{n}_v)$. 

\item
 $\sigma$ acts by sending $(v, e,
\vec{n_v})$ (resp. $(v, e, -\vec{n}_v)$)  to the dart which is
counter clockwise (resp. clockwise) from $(v, e)$ at $v$, maintaining
the orientation $\vec{n}_v$ (resp. $-\vec{n}_v$).

\item
$\tau$ acts by sending $(v, e, \pm \vec{n}_v)$ to 
$(\tilde{v}, e, \pm \vec{n}_v)$ when
$e$ is orientation preserving or
$(\tilde{v}, e, \mp \vec{n}_v)$ when $e$ is orientation reversing. 
\end{itemize}

On $\Omega_2$ we define the action of 
 $(\phi, \sigma, \tau)$ by:
\begin{itemize}
\item
$\phi$ acts by sending $(v, e, \pm \vec{n}_v; \pm)$ to 
$(v, e, \mp \vec{n}_v; \pm)$.

\item
$\sigma$ acts by sending $(v, e, \vec{n}_v; \pm)$ (resp. $(v, e,
  -\vec{n}_v; \pm)$) to the dart which is counter clockwise
  (resp. clockwise) from $(v, e, \pm)$ at $v$, maintaining the
  orientation $\vec{n}_v$ (resp. $-\vec{n}_v$).  

\item 
$\tau$ acts by sending $(v, e, \vec{n}_v; \pm)$ to $(v, e,
 \pm \vec{n}_v; \mp)$ 
(resp. sending $(v, e,
  -\vec{n}_v; \pm)$ to $(v, e, \mp \vec{n}_v; \mp)$ ) 
depending on whether $\tau$ is
  orientation preserving or reversing 
\end{itemize}

From this definition the following properties are apparent:
\begin{enumerate} 

\item $ \phi \sigma \phi = \sigma^{-1} $ and $\phi \tau \phi = \tau$

\item $\phi$ and $\tau$ are fixed point free permutations, each is a
  disjoint product of 2-cycles.  No 2-cycle of $\tau$ is a 2-cycle of $\phi$.  

\item $\sigma$ is a disjoint product of $2j_3$ 3-cycles, $2j_4$ 4-cycles,
  $\dots$, $2j_{2\nu}$ $2\nu$-cycles.

\item A connected map has the property that $\langle \phi, \sigma,
  \tau \rangle$ has a single orbit.

\item $\tau$ contains the data of how the darts of $K(D)$ are
  attached together to form the map; together with $\phi$, $\tau$
  determines how the local orientations relate globally. 

\item The cycles of $\sigma \circ \tau$  correspond to faces of
  the map.  There are two cycles of equal length
  for each face, representing the
  two directions the face can be traced along its edges.  

\end{enumerate}

This leads to the theorem:
\begin{thm}
There is a one-to-one correspondence between unoriented maps
 and triples $(\phi,
\sigma, \tau)$ satisfying conditions (1-6).  
\end{thm}

We will complete the proof by showing that given three permutations
$(\phi, \sigma, \tau)$ satisfying conditions (1-6), we can reconstruct
(up to the choice of local orientation) the map $( K(D), [\iota])$ and
the surface $X$.  

\begin{enumerate}[A.]

\item Let $K_0(D)$ be the set of cycle pairs $(\sigma_1, \sigma_2)$ in 
$\sigma$ such that  $\phi
  \sigma_1 \phi = \sigma_2^{-1}$.  

\item From each pair $(\sigma_1, \sigma_2)$ choose $\sigma_1$ to
  represent the counter clockwise or up orientation.

\item Let $K_1(D)$ be the set of 2-cycle pairs $(\tau_1, \tau_2)$ in
  $\tau$ such that $\phi \tau_1 \phi = \tau_2$.  
Each pair $(\tau_1, \tau_2)$ connects two (possibly the same) elements
  of $K_0(D)$.  If $(\tau_1, \tau_2)$ connects the up orientation to
  the up orientation this edge is orientation preserving, otherwise it
  is orientation reversing.  

\item To construct the surface $X$ and the embedding class $[\iota]$
  do the following:
\begin{enumerate}[(i)]
\item For each vertex $(\sigma_1, \sigma_2)$
  draw a polygon with sides labeled and oriented by the elements
  that $\sigma_1$ acts on.  

\item For an orientation preserving edge connecting two vertices, glue
  together the corresponding sides of their polygons so that the
  relative orientations agree at the corresponding side.  For an orientation
  reversing edge connecting two vertices, glue together the
  corresponding sides of their polygons so that the relative
  orientations are opposite at the side.  
\end{enumerate}

\end{enumerate}

In conclusion:  up to the choice of local orientations at each vertex, 
the data given
by $(\phi, \sigma, \tau)$ produces a map.  Likewise, a map, up to the
choice of local orientations at each vertex, produces a triple $(\phi,
\sigma, \tau)$ satisfying conditions 1-6.

\subsection{The algorithm}

We now have all the tools necessary to lift our algorithm to
the enumeration of unoriented maps.

Fix $\phi$ to be a fixed point free disjoint product of 2-cycles.  
Fix $\sigma$ to be a product of $2j_3$ 3-cycles, $2j_4$ 4-cycles,
$\dots$, and $2j_{2\nu}$ $2\nu$-cycles, such that $\phi$ connects
cycles in $\sigma$ into pairs $(\sigma_1, \sigma_2)$ of the same cycle
type, choose
$\sigma_1$ to represent the $\vec{n}_{(\sigma_1, \sigma_2)}$ orientation.

Cycle through all fixed point free products of disjoint 2-cycles $\tau$ such
that conditions 1-6 are satisfied by $(\phi, \sigma, \tau)$, in
particular no 2-cycle in $\phi$ appears in $\tau$.  If, for this
choice of $\tau$, the corresponding map is connected, then compute the
Euler characteristic:  the number of vertices is half the number of
cycles in $\sigma$, the number of edges is half the number of 2-cycles
in $\tau$, and the number of faces is half the number of cycles in
$\sigma \cdot \tau$.  Bin the results sorted by the Euler characteristic.

\section{Conclusion}

The algorithm presented here for counting oriented maps was used in
\cite{EMP} and the authors dissertation to evaluate constants of
integration appearing in the explicit calculation of $e_g(t_{2\nu})$
for low genus ($g \leq 3$).   Additionally the algorithm provided a
useful method for verifying the validity of the results in
\cite{EMP}.  

In looking beyond the calculations carried out in \cite{EMP} it became
apparent that it would be useful to have a similar algorithm to
compute the number of unoriented maps.  
The permutation structure associated to oriented maps has a
precise analogue in unoriented maps.  
This method may prove useful in further studies of the GOE and GSE
matrix integrals, and associated problems.

\end{document}